\documentclass[12pt,a4paper]{article}

\usepackage{amsmath,amsthm}
\usepackage{amssymb}
\usepackage{amsfonts}
\usepackage{amscd}

\newtheorem{thm}{Theorem}[section]
\newtheorem{theorem}[thm]{Theorem}

\newtheorem{corollary}[thm]{Corollary}
\newtheorem{lemma}[thm]{Lemma}

\newtheorem{proposition}[thm]{Proposition}

\theoremstyle{definition}

\newtheorem{definition}[thm]{Definition}

\newtheorem{remark}[thm]{Remark}
\newtheorem{observation}[thm]{Observation}

\newcommand{\cO}{\mathcal{O}}
\newcommand{\D}{\mathcal{D}}
\newcommand{\Det}{{D_q}}
\newcommand{\fg}{\mathfrak{g}}
\newcommand{\fgl}{{\mathfrak{gl}}}
\newcommand{\fsl}{{\mathfrak{sl}}}
\newcommand{\fp}{{\mathfrak{p}}}
\newcommand{\lra}{\longrightarrow}

\newcommand{\m}{{E}}
\newcommand{\e}{{\hbox{E}}}
\newcommand{\E}{{\hbox{e}}}
\newcommand{\f}{{\hbox{f}}}
\newcommand{\g}{{\hbox{g}}}

\newcommand{\Proof}{{\sl Proof}}
\newcommand{\Dij}{D^{1 \, \dots \, \widehat{j} \, \dots \, r \, i}}
\newcommand{\tdelta}{{\widetilde{\Delta}}}

\begin{document}

%
% Revised version of October 26th, 2008.
%
% \vskip 2cm
%

{\ }
\vskip-1cm

   \centerline{\LARGE \bf Quantum Duality Principle}
\bigskip
   \centerline {\LARGE\bf for Quantum Grassmannians}

\vskip 1cm

\centerline{R. Fioresi$^{\,\flat,}${\footnote{Partially supported by
the University of Bologna, funds for selected research topics.}},  \;
F. Gavarini$^{\,\sharp}$ }

\bigskip

\centerline{\it $^\flat$ Dipartimento di Matematica, Universit\`a
di Bologna }
\centerline{\it piazza di Porta S. Donato, 5  ---
I-40127 Bologna, ITALY}
\centerline{{\footnotesize e-mail: fioresi@dm.unibo.it}}

\bigskip

\centerline{\it $^\sharp$ Dipartimento di Matematica,
Universit\`a di Roma ``Tor Vergata'' }
\centerline{\it via della ricerca scientifica 1  ---
I-00133 Roma, ITALY}
\centerline{{\footnotesize e-mail: gavarini@mat.uniroma2.it}}

\vskip1,3cm

\begin{abstract}
The quantum duality principle (QDP) for homogeneous spaces gives four
recipes to obtain, from a quantum homogeneous space, a dual one, in the
sense of Poisson duality.  One of these recipes fails (for lack of the
initial ingredient) when the homogeneous space we start from is not a
quasi-affine variety.  In this work we solve this problem for the
quantum Grassmannian, a key example of quantum projective homogeneous
space, providing a suitable analogue of the QDP recipe.
\end{abstract}

\vskip21pt

\noindent   {\footnotesize {\bf Keywords:} \ {\sl Quantum Grassmann
Varieties}.}

\noindent   {\footnotesize 2000 {\it MSC:} \ Primary 20G42, 14M15;
Secondary 17B37, 17B62.}

\vskip31pt

\section{Introduction}
\label{intro}

   In the theory of quantum groups, the geometrical objects
that one takes into consideration are affine algebraic Poisson
groups and their infinitesimal counterparts, namely Lie bialgebras.
By ``quantization'' of either of these, one means a suitable
one-parameter deformation of one of the Hopf algebras associated
with them. They are respectively the algebra of regular function
$ \cO(G) \, $,  for a Poisson group  $ G $,  and the universal
enveloping algebra  $ U(\fg) $,  for a Lie bialgebra  $ \fg \, $.
Deformations of  $ \cO(G) $  are called {\sl quantum function
algebras}  (QFA), and are often denoted with  $ \cO_q(G) \, $,
while deformations of  $ U(\fg) $  are called {\sl quantum
universal enveloping algebras} (QUEA), denoted with
$ U_q(\fg) \, $.

\medskip

   The quantum duality principle (QDP), after its formulation in
\cite{ga2,ga3,ga4},  provides a recipe to get a QFA out of a QUEA,
and vice-versa. This involves a change of the underlying geometric
object, according to Poisson duality, in the following sense.
Starting from a QUEA over a Lie bialgebra  $ \, \fg = \text{\it
Lie}\,(G) \, $,  one gets a QFA for a dual Poisson group  $ G^* \, $.
Starting instead from a QFA over a Poisson group  $ G \, $,  one
gets a QUEA over the dual Lie bialgebra  $ \fg^* $.

\medskip

   In  \cite{cg},  this principle is extended to the wider
context of homogeneous Poisson  $ G $--spaces.  One describes
these spaces, in global or in infinitesimal terms, using
suitable subsets of  $ \cO(G) $  or of  $ U(\fg) \, $.  Indeed,
each homogeneous  $ G $--space  $ M $  can be realized as
$ G \big/ K $  for some closed subgroup  $ K $  of  $ G $
(this amounts to fixing a point in  $ M \, $:  it is shown in
\cite{cg},  \S 1.2, how to select such a point).  Thus we can
deal with either the space or the subgroup.  Now,  $ K $  can
be coded in infinitesimal terms by  $ U(\mathfrak{k}) $,  where
$ \, \mathfrak{k} := \text{\it Lie}\,(K) \, $,  and in global terms
by  $ \, \mathcal{I}(K) := \big\{\, \varphi \!\in\! \cO(G) \,\big|\,
\varphi(K) = 0 \big\} \, $,  the defining ideal of  $ K \, $.  Instead,
$ G \big/ K $  can be encoded infinitesimally by  $ \, U(\fg) \,
\mathfrak{k} \, $ and globally by  $ \, \cO\big(G\big/K\big) \equiv
{\cO(G)}^K $,  the algebra of  $ K $--invariants in  $ \cO(G) \, $.
Note that  $ \, U(\fg) \big/ U(\fg) \, \mathfrak{k} \, $  identifies
with the set of left-invariant differential operators on  $ G\big/K \, $,
or the set of  $ K $--invariant,  left-invariant differential operators
on  $ G \, $.
                                        \par
   These constructions  {\sl all\/}  make sense in formal
geometry, i.e.~when dealing simply with formal groups and formal
homogeneous spaces, as in  \cite{cg}.  Instead, if one looks for
{\sl global\/}  geometry, then one construction might fail, namely
the description of  $ G \big/ K \, $  via its function algebra
$ \, \cO\big(G\big/K\big) = {\cO(G)}^K \, $.  In fact, this makes
sense   --- i.e.,  $ {\cO(G)}^K $  is enough to describe  $ G \big/ K $
---   if and only if the variety  $ G \big/ K $  is  {\sl quasi-affine}.
In particular, this is not the case if  $ G \big/ K $  is projective,
like, for instance, when  $ G \big/ K $  is a Grassmann variety.
                                        \par
   By ``quantization'' of the homogeneous space  $ G \big/ K $  one
means any quantum deformation (in suitable sense) of any one of the
four algebraic objects mentioned before which describe either
$ G \big/ K $  or  $ K \, $.  Moreover one requires that given
an infinitesimal or a global quantization for the group  $ G $,
denoted by  $ U_q(\fg) $  or  $ \cO_q(G) $  respectively, the
quantization of the homogeneous space admits a  $ U_q(\fg) $--action
or a  $ \cO_q(G) $--coaction  respectively, which yields a quantum
deformation of the algebraic counterpart of the  $ G $--action
on  $ G \big/ K \, $.
                                        \par
   The QDP for homogeneous  $ G $--spaces  (cf.~\cite{cg})  starts
from an infinitesimal (global) quantization of a  $ G $--space,  say
$ G \big/ K $,  and provides a global (infinitesimal) quantization
for the Poisson dual  $ G^* $--space.  The latter is  $ \, G^* \!
\big/ K^\perp $  (with  $ \, \text{\it Lie}\,\big(K^\perp\big)
= \mathfrak{k}^\perp \, $,  the orthogonal subspace   --- with
respect to the natural pairing between  $ \fg $  and its dual space
$ \fg^* $  ---   to  $ \mathfrak{k} $  inside  $ \fg^* \, $).  In
particular, the principle gives a concrete recipe
  $$  \cO_q\big(G\big/K\big) \;\;
\circ\hskip-5,3pt\relbar\relbar \relbar \relbar\joinrel\rightsquigarrow
\;\;  {\cO_q\big(G\big/K\big)}^\vee =: U_q\big(\mathfrak{k}^\perp\big)  $$
in which the right-hand side is a quantization of  $ \, U \big(
\mathfrak{k}^\perp \big) \, $.
                                        \par
   However, this recipe makes no sense when  $ \cO_q\big(G\big/K\big) $
is not available.  In the non-formal setting this is the case whenever
$ G \big/ K $  is not quasi-affine, e.g.~when it is projective.

\medskip

   In this paper we show how to solve this problem in the special case
of the Grassmann varieties, taking  $ G $  as the general linear group
and  $ \, K = P \, $  a maximal parabolic subgroup.  We adapt the basic
ideas of the original QDP recipe to these new ingredients, and we obtain
a new recipe
  $$  \cO_q \big( G \big/ P \big)  \;\;
\circ\hskip-5,3pt\relbar \relbar \relbar \relbar\joinrel\rightsquigarrow
\;\;  \widehat{{\cO_q \big( G \big/ P \big)}^\vee}  $$
which perfectly makes sense, and yields the same kind of result as
predicted by the QDP for the quasi-affine case.  In particular,
$ \widehat{{\cO_q \big( G \big/ P \big)}^\vee} $  is a quantization
of  $ U\big(\mathfrak{p}^\perp\big) \, $,  obtained through a
$ (q-1) $--adic  completion  process.

\medskip

   Our construction goes as follows.

\medskip

   First, we consider the embedding of  the Grassmannian $ G \big/ P $
(where  $ \,  G := {GL}_n \, $  or  $ \, G := {SL}_n \, $, and  $ P $
is a parabolic subgroup of  $ G \, $)  inside a projective space, given
by Pl{\"u}cker coordinates.  This will give us the first new ingredient:
  $$  \cO\big( G \big/ P \big) \; := \; \hbox{ring of homogeneous
coordinates on} \  G \big/ P  \quad .  $$
   \indent   Many quantizations  $ \cO_q \big( G \big/ P \big) $  of
$ \cO\big( G \big/ P \big) $  already exist in the literature (see, e.g.,
\cite{fi1, lr, tt}).  All these quantizations, which are equivalent,
come together with a quantization of the natural  $ G $--action
on  $ \, G/P \, $.

\medskip

   In the original recipe  (see \cite{cg})  $ \; \cO_q \big( G\big/K \big)
\, \circ\hskip-5,3pt\relbar \relbar \relbar \relbar\joinrel\rightsquigarrow
{\cO_q\big( G\big/K \big)}^\vee \; $  of the QDP (when  $ G \big/ K $  is
quasi affine) we need to look at a neighborhood of the special
point  $ eK $  (where  $ \, e \in G \, $  is the identity), and at a
quantization of it.  Therefore, we shall replace the projective variety
$ G \big/ P $  with such an affine neighborhood, namely the big cell of
$ G \big/ P \, $.  This amounts to realize the algebra of regular functions
on the big cell as a ``homogeneous localization'' of  $ \cO\big( G\big/P
\big) $,  say  $ \cO^{\,\text{\it loc}} \big( G \big/ P \big) $,  by
inverting a suitable element.  We then do the same at the quantum level,
via the inversion of a suitable almost central element in  $ \cO_q\big(
G\big/P \big) $   --- which lifts the previous one in  $ \cO\big( G\big/P
\big) \, $.  The result is a quantization  $ \cO_q^{\,\text{\it loc}}
\big( G \big/ P \big) $  of the coordinate ring of the big cell.

\medskip

   Hence we are able to  {\sl define\/}  $ \, {\cO_q\big( G \big/ P
\big)}^\vee \! := {\cO_q^{\,\text{\it loc}}\big( G \big/ P \big)}^\vee $,
where the right-hand side is given by the original QDP recipe applied to
the big cell as an affine variety (we can forget any group action at
this step).  By the very construction, this  $ {\cO_q \big( G \big/ P
\big)}^\vee $  should be a quantization of  $ U \big( \mathfrak{p}^\perp
\big) \, $  (as an algebra).  Indeed, we prove that this is the case, so we
might think at  $ {\cO_q \big( G \big/ P \big)}^\vee $  as a quantization
(of infinitesimal type) of the variety  $ G^* \big/ P^\perp \, $.  On the
other hand, the construction does not ensure that  $ {\cO_q\big( G \big/
P \big)}^\vee $  also admits a quantization of the  $ G^* $--action
on  $ G^* \big/ P^\perp \, $  (just like the big cell is not a
$ G $--space).  As a last step, we look at  $ \widehat{{\cO_q \big(
G \big/ P \big)}^\vee} $,  the  $ (q \! - \! 1) $--adic  completion
of  $ {\cO_q\big( G \big/ P \big)}^\vee $.  Of course, it is again
a quantization of  $ U \big( \mathfrak{p}^\perp \big) \, $  (as an
algebra).  But in addition, it admits a coaction of the  $ (q \! -
\! 1) $--adic  completion of  $ \cO_q(G)^\vee $   --- which is a
quantization of  $ U(\fg^*) $. This coaction  yields a quantization
of the infinitesimal  $ G^* $--action  on  $ G^* \big/ P^\perp $.
Therefore, in a nutshell,  $ \, \widehat{{\cO_q\big( G \big/ P \big)}^\vee}
\, $  is a quantization of  $ \, G^* \big/ P^\perp \, $  {\sl as a
homogeneous  $ G^* $--space},  in the sense explained above.

\medskip

   Notice that our arguments could be applied to any  {\sl projective\/}
homogeneous  $ G $--space  $ X \, $,  {\sl up to having the initial data
to start with}.  Namely, one needs an embedding of  $ X $  inside a
projective space, a quantization (compatible with the  $ G $--action)
of the ring of homogeneous coordinates of  $ X $  ({w.r.t.} such an
embedding), and a quantization of a suitable open dense affine subset
of  $ X \, $.  This program is carried out in detail in a separate
work (see  \cite{cfg}).

\medskip

   Finally, this paper is organized as follows.
                                     \par
   In section  \ref{quantum}  we fix the notation, and we describe the
Manin deformations of the general linear group (as a Poisson group), and
of its Lie bialgebra, together with its dual.  In section  \ref{grass}
we briefly recall results concerning the constructions of the quantum
Grassmannian  $ \cO_q \big( G \big/ P \big) $  and its quantum big cell
$ \cO_q^{\,\text{\it loc}} \big( G \big/ P \big) \, $.  These are known
results, treated in detail in  \cite{fi1,fi2}.  Finally, in section
\ref{qdpgrass}  we extend the original QDP to build  $ {\cO_q \big(
G \big/ P \big)}^\vee $,  and we show that its  $ (q - 1) $--adic
completion is a quantization of the homogeneous  $ G^* $--space
$ G^* \big/ P^\perp $  dual to the Grassmannian  $ G \big/ P \, $.

\bigskip
\medskip

   \centerline{\bf Acknowledgements}

\medskip

   The first author wishes to thank the Dipartimento di Matematica ``Tor
Vergata'', and in particular Prof.~V.~Baldoni and Prof.~E.~Strickland, for
the warm hospitality during the period in which this paper was written.
                                                    \par
   Both authors also thank D.~Parashar and M.~Marcolli for their kind
invitation to the workshop ``Quantum Groups and Noncommutative
Geometry'' held at MPIM in Bonn during August 6--8, 2007.

\bigskip
\bigskip

\section{The Poisson Lie group  $ {GL}_n(\Bbbk) $  and its quantum
deformation}
\label{quantum}

\bigskip

   {} \indent   Let  $ \Bbbk $  be any
%
% algebraically closed
%
 field of characteristic zero.

\medskip

   In this section we want to recall the construction of a quantum
deformation of the Poisson Lie group  $ \, {GL}_n := {GL}_n(\Bbbk) \, $.
We will also describe explicitly the bialgebra structure of its Lie
algebra  $ \, \fgl_n := \fgl_n(\Bbbk) \, $  in a way that fits our
purposes, that is to obtain a quantum duality principle for the
Grassmann varieties for  $ {GL}_n $  (see \S  \ref{qdpgrass}).

\medskip

   Let  $ \, \Bbbk_q = \Bbbk\big[q,q^{-1}\big] \, $  (where  $ q $  is
an indeterminate), the ring of Laurent polynomials over  $ q \, $,  and
let  $ \, \Bbbk(q) \, $  be the field of rational functions in  $ q \, $.

\medskip

\begin{definition} \label{manin}
   The {\sl quantum matrix algebra} is defined as
  $$  \cO_q(M_{m \times n}) \; = \; \Bbbk_q \big\langle {\{\, x_{ij}
\,\}}_{1 \leq i \leq m}^{1 \leq j \leq n} \,\big\rangle \Big/ I_M  $$
where the  $ x_{ij} $'s  are non commutative indeterminates, and
$ I_M $  is the two-sided ideal generated by the {\sl Manin relations}
  $$  \displaylines{
   x_{ij} \, x_{ik} \; = \; q \, x_{ik} \, x_{ij} \; ,  \qquad
 x_{ji} \, x_{ki} \; = \; q \, x_{ki} \, x_{ji}  \; \qquad
\forall \;\; j < k  \cr
   {\ }  \qquad \qquad
   x_{ij} \, x_{kl} \; = \; x_{kl} \, x_{ij}  \; \qquad \qquad  \forall
\;\; i < k \, , \, j > l  \hbox{\ \ or \ }  i > k \, , \, j < l  \cr
   x_{ij} \, x_{kl} \, - \, x_{kl} \, x_{ij} \; = \; \big( q - q^{-1} \big)
\, x_{kj} \, x_{il}  \; \qquad  \forall \;\; i < k \, , \, j < l
\cr }  $$
\end{definition}

\medskip

   {\bf Warning:}  sometimes these relations appear with  $ q $
exchanged with  $ q^{-1} \, $.

\bigskip

   For simplicity we will denote  $ \cO_q(M_{n \times n}) $  with
$ \cO_q(M_{n}) \, $.

\medskip

   There is a coalgebra structure on  $ \cO_q(M_{n}) \, $,
given by
  $$  \Delta(x_{ij}) \; = \; {\textstyle \sum\limits_{k=1}^n} \, x_{ik}
\otimes x_{kj} \quad ,  \qquad  \epsilon(x_{ij}) \; = \; \delta_{ij}
\qquad   \eqno ( \, 1 \leq i \, , j \leq n \,)  $$

\medskip

   The  {\sl quantum general linear group\/}  and the  {\sl quantum
special linear group\/}  are defined in the following way:
  $$  \cO_q(GL_n) \! := \cO_q(M_n)[T] \Big/ \! \big( T \Det - 1 \, ,
1 - T \Det \big) \; ,  \;\;  \cO_q(SL_n) \! := \cO_q(M_n) \Big/ \!
\big( \Det - 1 \big)  $$
where  $ \; \Det := \sum_{\sigma \in \mathcal{S}_n} \, (-q)^{\ell(\sigma)}
\, x_{1 \, \sigma(1)} \cdots x_{n \, \sigma(n)} \; $  is a central
element, called the {\sl quantum determinant}.

\bigskip

   {\bf Note:}  We use the same letter to denote the generators
$ x_{ij} $  of  $ \cO_q(M_{m \times n}) \, $,  of  $ \cO_q(GL_n) $
and of  $ \cO_q(SL_n) \, $:  the context will make clear where they
sit.

\bigskip

   The algebra  $ \cO_q({GL}_n) $  is a quantization of the algebra
$ \cO({GL}_n) $  of regular functions on the affine algebraic group
$ {GL}_n \, $,  in the following sense:  $ \; \cO_q({GL}_n) \big/
(q \! - \! 1) \, \cO_q({GL}_n) \; $  is isomorphic to  $ \cO({GL}_n) $
as a Hopf algebra (over the field $ \Bbbk \, $).  Similarly, $ \cO_q
({SL}_n) $  is a quantization of the algebra  $ \cO({SL}_n) $  of
regular functions on $ {SL}_n \, $.  Both  $ \cO_q({GL}_n) $  and
$ \cO_q({SL}_n) $  are Hopf algebras, that is, they also have the
antipode.  For more details on these constructions see for example
\cite{cp},  pg. 215.

\medskip

   By general theory,  $ \cO({GL}_n) $  inherits from
$ \cO_q({GL}_n) $  a Poisson bracket, which makes it into a
Poisson Hopf algebra, so that  $ {GL}_n $  becomes a Poisson
group.  We want to describe now its Poisson bracket.  Recall that
  $$  \cO({GL}_n) \; = \; \Bbbk \big[ {\{\, \bar{x}_{ij} \,\}}_{i, j
= 1, \dots, n} \,\big][t] \Big/ \big( t \, d - 1 \big)  $$
where  $ \; d := \text{\sl det}\, \big( \bar{x}_{i,j} \big)_{i, j =
1, \dots, n} \; $  is the usual determinant.  Setting  $ \, \bar{x}
= \pi(x) \, $  for  $ \; \pi : \cO_q(GL_n) \lra \cO(GL_n) \; $,  \,
the Poisson structure is given (as usual) by
  $$  \big\{ \bar{a} \, , \bar{b} \,\big\} \; := \; {(q-1)}^{-1}
\, ( a \, b - b \, a )\Big|_{q=1}   \eqno \forall \;\; \bar{a} \, ,
\bar{b} \in \cO({GL}_n) \;\; .  $$
In terms of generators, we have
  $$  \displaylines{
   \big\{ \bar{x}_{ij} \, , \bar{x}_{ik} \big\} \, = \, \bar{x}_{ij}
\, \bar{x}_{ik}  \quad  \forall \;\; j < k \; ,   \phantom{\Big|}
 \quad  \big\{ \bar{x}_{ij} \, , \bar{x}_{\ell k} \big\} \, =
\, 0  \qquad  \forall \;\; i < \ell \, , k < j  \cr
   \big\{ \bar{x}_{ij} \, , \bar{x}_{\ell j} \big\} \, = \,
\bar{x}_{ij} \, \bar{x}_{\ell j}  \quad  \forall \;\; i < \ell \, ,
\phantom{\Big|}  \qquad  \big\{ \bar{x}_{ij} \, , \bar{x}_{\ell k} \big\}
\, = \, 2 \, \bar{x}_{ij} \, \bar{x}_{\ell k}  \quad  \forall \;\;
i < \ell \, , j < k  \cr
   \big\{ d^{-1} , \bar{x}_{ij} \big\} \, = \,0 \; ,
\phantom{\Big|}  \quad  \big\{ d \, , \bar{x}_{ij} \big\}
\,= \, 0  \qquad \forall \;\; i, j = 1, \dots, n \, .  \cr }  $$
   \indent   As  $ {GL}_n $  is a Poisson Lie group, its Lie algebra
$ \fgl_n $  has a Lie bialgebra structure (see  \cite{cp},  {pg.~\!}24).
To describe it, let us denote with  $ \m_{ij} $  the elementary
matrices, which form a basis of  $ \fgl_n \, $.
%
%  $$  \m_{i,j} := {\big( \delta_{\ell,i} \, \delta_{j,s}
% \big)}_{\ell, = 1, \dots, n}^{s = 1, \dots, n} \, \qquad
% \forall \; i, j = 1, \dots, n  $$
%
Define  ($ \, \forall \; i = 1, \dots, n-1 \, $,  $ \, j = 1, \dots,
n \, $)
  $$  e_i := \m_{i,i+1} \; ,  \quad  g_j := \m_{j,j} \; ,  \quad
f_i := \m_{i+1,i} \; ,  \quad  h_i := g_{i} - g_{i+1}  $$
%
% (for all  $ \; i = 1, \dots, n-1 \, $,  $ \, j = 1, \dots, n \, $).
%
 Then
$ \; \big\{\, e_i \, , \, f_i \, , \, g_j \;\big|\; i = 1, \dots, n-1,
\, j = 1, \dots, n \,\big\} \; $  is a set of Lie algebra generators of
$ \fgl_n \, $,  and a Lie cobracket is defined on  $ \fgl_n $  by
%
%   $$  \delta(e_i) \, = \, h_i \otimes e_i - e_i \otimes h_i \, ,  \quad
% \delta(g_j) \, = \, 0 \, ,  \quad  \delta(f_i) \, = \, h_i \otimes f_i
% - f_i \otimes h_i \;   \eqno \forall \;\, i , j .  $$
% %
%
  $$  \delta(e_i) \, = \, h_i \wedge e_i \;\; ,  \quad
\delta(g_j) \, = \, 0 \;\; ,  \quad  \delta(f_i) \, = \, h_i \wedge f_i
\;   \eqno \forall \;\, i , j .  $$
This cobracket makes  $ \fgl_n $  itself into a  {\sl Lie bialgebra\/}:
this is the so-called  {\sl standard\/}  Lie bialgebra structure on
$ \fgl_n \, $.  It follows immediately that  $ \, U(\fgl_n) \, $  is
a co-Poisson Hopf algebra, whose co-Poisson bracket is the (unique)
extension of the Lie cobracket of  $ \fgl_n $  while the Hopf
structure is the standard one.

\medskip

   Similar constructions hold for the group  $ {SL}_n \, $.  One simply
drops the generator  $ d^{-1} \, $,  imposes the relation  $ \, d \!
= \! 1 \, $,  in the description of  $ \cO(SL_n) \, $,  and replaces the
$ g_s $'s  with the  $ h_i $'s  ($ \, i = 1, \dots, n \, $)  when
describing  $ \, \fsl_n \, $.

\medskip

   Since  $ \fgl_n $ is a Lie bialgebra, its dual space  $ \fgl_n^{\,*} $
admits a Lie bialgebra structure, dual to the one of  $ \fgl_n \, $.  Let
$ \, \big\{\, \e_{ij} := \m_{ij}^{\,*} \;\big|\; i, j = 1, \dots, n
\,\big\} \, $  be the basis of  $ \, \fgl_n^{\,*} \, $  dual to the
basis of elementary matrices for  $ \fgl_n \, $.  As a Lie algebra,
$ \fgl_n^{\,*} $  can be realized as the subset of  $ \, \fgl_n
\oplus \fgl_n \, $  of all pairs
  $$  \left( \! \begin{pmatrix}
   \! -m_{11} \!  &  \! 0 \!  &  \! \cdots \!  & \!   0 \!   \\
   \! m_{21} \!  &  \! -m_{22} \!  &  \! \cdots & \!  0 \!   \\
   \! \vdots \!  &  \! \vdots \!  &  \! \vdots \!  &  \! \vdots \!   \\
   \! m_{n-1,1} \!  &  \! m_{n-1,2} \!  &  \! \cdots \!  &  \! 0 \!   \\
   \! m_{n,1} \!  &  \! m_{n,2} \!  &  \! \cdots \!  &  \! -m_{n,n}
      \end{pmatrix} , \,
      \begin{pmatrix}
   m_{11} \!  &  \! m_{12} \!  &  \! \cdots \!  &  \! m_{1,n-1} \!
&  \! m_{1,n} \!   \\
   \! 0 \!  &  \! m_{22} \!  &  \! \cdots \!  &  \! m_{2,n-1} \!  &
\! m_{2,n} \!   \\
   \! \vdots \!  &  \! \vdots \!  &  \! \vdots \!  &  \! \vdots \!  &
\! \vdots \!   \\
   \! 0 \!  &  \! 0 \!  &  \! \cdots \!  &  \! m_{n-1,n-1} \!  &
\! m_{n-1,n} \!   \\
   \! 0 \!  &  \! 0 \!  &  \! \cdots \!  &  \! 0 \!  &  \! m_{n,n} \!
      \end{pmatrix} \! \right)  $$
with its natural structure of Lie subalgebra of  $ \, \fgl_n \oplus
\fgl_n \, $.  In fact, the elements  $ \, \e_{ij} \, $  correspond to
elements in  $ \, \fgl_n \oplus \fgl_n \, $  in the following way:
  $$  \e_{ij} \cong \big( \m_{ij} \, , 0 \big)  \hskip5pt  \forall \;
i \! > \! j \, ,  \hskip7pt  \e_{ij} \cong \big(\! -\m_{ij} \, , +\m_{ij}
\big)  \hskip5pt  \forall \; i \! = \! j \, ,  \hskip7pt  \e_{ij} \cong
\big( 0 \, , \m_{ij} \big)  \hskip5pt  \forall \; i < j \, .  $$
Then the Lie bracket of  $ \, \fgl_n^{\,*} \, $  is given by
  $$  \begin{array}{clc}
   \big[ \e_{i,j} \, , \, \e_{h,k} \big]  \hskip-5pt  &  = \, \delta_{j,h}
\, \e_{i,k} - \delta_{k,i} \, \e_{h,j} \;\; ,   &   \forall  \;\; i \!
\leq \! j \, , \, h \! \leq \! k  \;\;\,  \text{and}  \;\; \forall \;\;
i \! > \! j \, , \, h \! > \! k  \\
                                 \\
   \big[ \e_{i,j} \, , \, \e_{h,k} \big]  \hskip-5pt  &  = \, \delta_{k,i}
\, \e_{h,j} - \delta_{j,h} \, \e_{i,k} \;\; ,  &   \forall \;\; i \! = \!
j \, , \, h \! > \! k  \;\;\,  \text{and}  \;\; \forall \;\; i \! > \! j
\, , \, h \! = \! k  \\
                                 \\
   \big[ \e_{i,j} \, , \, \e_{h,k} \big]  \hskip-5pt  &  = \, 0 \; ,  &
\forall\;\; i \! < \! j \, , \, h \! > \! k  \;\;\,  \text{and}  \;\;
\forall \;\; i \! > \! j \, , \, h \! < \! k
      \end{array}  $$
   \indent   Note that the elements  ($ \, 1 \leq i \leq n \! - \! 1 \, $,
$ \, 1 \leq j \leq n \, $)
  $$  \E_i \, = \, e_i^{\,*} \, = \, \e_{i,i+1} \;\; ,  \qquad
\f_i \, = \, f_i^{\,*} \, = \, \e_{i+1,i} \;\; ,  \qquad
\g_j \, = \, g_j^{\,*} \, = \, \e_{jj}  $$
%
% ($ \, 1 \! \leq \! i \! \leq \! n \! - \! 1 \, $,  $ \, 1 \! \leq \! j
% \! \leq \! n \, $)
%
 are Lie algebra generators of  $ \fgl_n^{\,*} \, $.
In terms of them, the Lie bracket reads
  $$  \big[ \E_i \, , \f_j \big] \, = \, 0 \; ,  \qquad  \big[ \g_i \, ,
\E_j \big] \, = \, \delta_{ij} \, \E_i \; ,  \qquad  \big[ \g_i \, , \f_j
\big] \, = \, \delta_{ij} \, \f_j   \eqno \forall \;\; i, j \; .  $$

\medskip

   On the other hand, the Lie cobracket structure of  $ \fgl_n^{\,*} $
is given by
  $$  \delta\big(\e_{i,j}\big) \, = \,
{\textstyle \sum\limits_{k=1}^n} \, \e_{i,k} \wedge \e_{k,j}
\eqno \forall \;\; i, j = 1, \dots, n  $$
where
%
% hereafter we use shorthand notation
%
 $ \, x \wedge y := x \otimes y - y \otimes x \; $.

\medskip

   Finally, all these formul{\ae}  also provide a presentation of
$ U\big(\fgl_n^{\,*}\big) $  as a co-Poisson Hopf algebra.

\medskip

   A similar description holds for  $ \, \fsl_n^{\,*} \! = \fgl_n^{\,*}
\Big/ Z\big(\fgl_n^{\,*}\big) \, $,  where  $ \, Z \big( \fgl_n^{\,*}
\big) \, $  is the centre of  $ \fgl_n^{\,*} \, $,  generated by  $ \,
\mathfrak{l}_n := \g_1 \! + \cdots + \g_n \, $.  The construction is
immediate by looking at  the embedding  $ \, \fsl_n \hookrightarrow
\fgl_n \, $.

\bigskip

\section{\hskip-9pt The quantum Grassmannian and its big cell}
\label{grass}

   {} \indent   In this section we want to briefly recall the
construction of a quantum deformation of the Grassmannian of
$ r $--spaces  inside an  $ n $--dimensional  vector space and
its big cell, as they appear in  \cite{fi1,fi2}.  The
quantum Grassmannian ring will be obtained as a quantum homogeneous
space, namely its deformation will come together with a deformation of
the natural coaction of the function algebra of the general linear
group on it.  The deformation will also depend on a specific embedding
(the Pl{\"u}cker one) of the Grassmann variety into a projective
space.  This deformation is very natural, in fact it embeds into
the deformation of its big cell ring.  Let's see explicitly these
constructions.

\medskip

   Let  $ \, G := {GL}_n \, $,  and let  $ P $  and  $ P_1 $  be the
standard parabolic subgroups
  $$  P \, := \, \left\{ \begin{pmatrix} A & B \\ 0 & C \end{pmatrix}
\right\} \subset {GL}_n \quad ,  \qquad  P_1 \, := \, P \, {\textstyle
\bigcap} \, {SL}_n  $$
where  $ A $  is a square matrix of size  $ r \, $,  with  $ \, 0 < r
< n \, $.

\medskip

\begin{definition}
  The  {\sl quantum Grassmannian coordinate ring}  $ \, \cO_q \big(
G \big/ P \big) \, $  with respect to the Pl{\"u}cker embedding is
the subalgebra of  $ \cO_q({GL}_n) $  generated by the quantum minors
(called  {\sl quantum Pl{\"u}cker coordinates\/})
  $$  D^I \, = \, D^{i_1 \dots i_r} \, := \, {\textstyle \sum\limits_{\sigma \in
\mathcal{S}_r}} \, {(-q)}^{\ell(\sigma)} \, x_{i_1 \, \sigma(1)}
\, x_{i_2 \, \sigma(2)} \cdots x_{i_r \, \sigma(r)}   \quad .  $$
for every ordered  $ r $--tuple  of indices  $ \; I =\{i_1 < \cdots < i_r\} \, $.
\end{definition}

\smallskip

\noindent
 {\it  $ \underline{\text{Remark}} $:}  \,  Equivalently,  $ \,
\cO_q \big( G\big/P \big) \, $  may be defined in the same way but
with  $ \cO_q({SL}_n) $  instead of  $ \cO_q({GL}_n) \, $.

\bigskip

   The algebra  $ \cO_q\big(G\big/P\big) $  is a quantization of the
Grassmannian  $ G\big/P $  in the usual sense: the  $ \Bbbk $--algebra
$ \, \cO_q\big( G\big/P\big) \Big/ (q-1) \, \cO_q\big(G\big/P\big)
\, $  is isomorphic to  $ \cO\big(G\big/P\big) \, $,  the algebra
of homogeneous coordinates of  $ G\big/P $  with respect to the
Pl{\"u}cker embedding.
%
% (thought of as a projective subvariety of a suitable projective space).
%
In addition,  $ \cO_q\big(G\big/P\big) $  has an important property
w.r.t.~$ \cO_q(G) \, $,  given by the following result:

\medskip

\begin{proposition}     \label{firstintersection}
  $$  \cO_q\big(G\big/P\big) \,\;{\textstyle \bigcap}\;\,
(q-1) \, \cO_q(G) \,\; = \;\, (q-1) \; \cO_q\big(G\big/P\big)  $$
\end{proposition}

\Proof.  By Theorem 3.5 in  \cite{tt},  we have that certain products of
minors  $ \, {\{p_i\}}_{i \in I} \, $  form a basis of  $ \cO_q\big( G
\big/ P \big) \, $  over  $ \Bbbk_q \, $.  Thus, a generic element in
$ \, \cO_q\big(G\big/P\big) \,\bigcap\, (q \! - \! 1) \, \cO_q(G) \, $
can be written as
  $$  {\textstyle \sum_{i \in I}} \, \alpha_i \, p_i \; = \;
(q-1) \, \phi   \eqno (3.1)  $$
for some  $ \, \phi \in \cO_q(G) \, $.  Moreover, the specialization
map
  $$  \pi_G \; \colon \; \cO_q(G) \;
\relbar\joinrel\relbar\joinrel\relbar\joinrel\twoheadrightarrow
\; \cO_q(G) \Big/ (q-1) \, \cO_q(G) \; = \; \cO(G)  $$
maps  $ \, {\{p_i\}}_{i \in I} \, $  onto a basis  $ \, {\big\{
\pi_G(p_i) \big\}}_{i \in I} \, $  of  $ \cO\big( G\big/P \big) \, $,
the latter being a subalgebra of  $ \cO(G) \, $.  Therefore, applying
$ \pi_G $  to (3.1) we get  $ \, \sum_{i \in I} \overline{\alpha_i} \,
\pi_G(p_i) = 0 \, $,  where  $ \, \overline{\alpha_i} := \alpha_i \mod
(q \! - \! 1) \, \Bbbk_q \, $,  for all  $ \, i \in I \, $.  This forces
$ \, \alpha_i \in (q \! - \! 1) \, \Bbbk_q \, $  for all  $ i \, $,  by
the linear independence of the  $ \pi_G(p_i) $'s,  whence the claim.
\qed

\bigskip

   An immediate consequence of Proposition  \ref{firstintersection}
is that the canonical map
  $$  \cO_q\big( G\big/P \big) \Big/ (q-1) \, \cO_q\big( G\big/P \big)
\; \relbar\joinrel\relbar\joinrel\relbar\joinrel\longrightarrow \;
\cO_q(G) \Big/ (q-1) \, \cO_q(G)  $$
is  {\sl injective}.  Therefore, the specialization map
  $$  \pi_{G/P} \; \colon \; \cO_q\big( G\big/P \big)
\;\relbar\joinrel\relbar\joinrel\relbar\joinrel\twoheadrightarrow\;
\cO_q\big( G\big/P \big) \Big/ (q-1) \, \cO_q\big( G\big/P \big)  $$
coincides with the restriction to  $ \cO_q\big( G\big/P \big) $  of
the specialization map
  $$  \pi_G \; \colon \; \cO_q(G) \;
\relbar\joinrel\relbar\joinrel\relbar\joinrel\twoheadrightarrow
\; \cO_q(G) \Big/ (q-1) \, \cO_q(G)  \quad .  $$
   \indent   Moreover   --- from a geometrical point of view ---   the
key consequence of this property is that  {\sl  $ P $  is a coisotropic
subgroup\/}  of the Poisson group  $ G \, $.  This implies the existence
of a well defined Poisson structure on the algebra  $ \, \cO\big( G\big/P
\big) \, $,  inherited from the one in  $ \cO(G) \, $.

\bigskip

\begin{observation} \label{comultiplication}
  The quantum deformation  $ \cO_q\big( G\big/P \big) $
                   comes naturally\break
   equipped with a coaction of  $ \cO_q({GL}_n) $   --- or, similarly,
of  $ \cO_q({SL}_n) $  ---   on it, obtained by restricting the
comultiplication  $ \Delta \, $.  This reads
  $$  \begin{array}{cccc}
   \Delta{\big|}_{\cO_q(G/P)} :  &  \cO_q\big( G \big/ P \big)  &  \lra
&  \cO_q(G) \otimes \cO_q\big( G \big/ P \big)  \\
   &  D^I  &  \mapsto  &  \sum_K D^I_K \otimes D^K
      \end{array}  $$
where, for any  $ \, I = (i_1 \dots i_r) \, $,  $ \, K = (k_1 \dots k_r)
\, $,  with  $ \, 1 \leq i_1 < \dots < i_r \leq n \, $,  $ \, 1 \leq k_1
< \dots < k_r \leq n \, $,  we denote by  $ D^I_K $  the  {\sl quantum
minor}
  $$  D^I_K \; \equiv \; D^{i_1 \dots i_r}_{k_1 \dots k_r}  \;
:= \;  {\textstyle \sum\limits_{\sigma \in \mathcal{S}_r}} \,
{(-q)}^{\ell(\sigma)} \, x_{i_1 \, k_{\sigma(1)}} \, x_{i_2 \,
k_{\sigma(2)}} \cdots x_{i_r \, k_{\sigma(r)}}  \quad .  $$
This provides a quantization of the natural coaction of
$ \cO(G) $  onto  $ \cO\big( G \big/ P \big) \, $.
\end{observation}

\medskip

   The ring  $ \cO_q\big( G\big/P \big) $  has been fully described in
\cite{fi1}  in terms of generators and relations.  We refer the reader
to this work for further details.

\medskip

   We now turn to the construction of the quantum big cell ring.

\medskip

\begin{definition}
  Let  $ \, I_0 = (1 \dots r) \, $,  $ \, D_0 := D^{I_0} \, $.  Define
  $$  \cO_q(G)\big[D_0^{-1}\big] \; := \;
\cO_q(G)[T] \Big/ \big( T \, D_0 - 1 \, , D_0 \, T - 1 \big)  $$
Moreover, we define the  {\sl big cell ring\/}  $ \, \cO_q^{\,\text{\it
loc}}\big( G\big/P \big) \, $  to be the  $ \Bbbk_q $--subalgebra  of
$ \, \cO_q(G)\big[D_0^{-1}\big] \, $  generated by the elements
  $$  {\ } \qquad  t_{ij} \, := \, {(-q)}^{r-j} \, \Dij \, D_0^{-1}
\eqno \forall \;\;\, i \, , j \; : \; 1 \leq j \leq r < i \leq n \; .  $$
%
% where  $ \, I = \big(1 \dots \widehat j \dots r \, i \,\big) \, $.
%
\end{definition}

See \cite{fi2} for more details.

\medskip

   As in the commutative setting, we have the following result:

\medskip

\begin{proposition} \label{classical}
  $ \; \displaystyle{ \cO_q^{\,\text{\it loc}}\big( G\big/P \big)
\, \cong \, \cO_q\big( G\big/P \big)\big[D_0^{-1}\big]_{proj} } \; $,
\; where the right-hand side is the degree-zero component of
$ \; \cO_q\big( G\big/P \big)[T] \Big/ \big(T D_0 - 1 \, ,
D_0 \, T - 1 \big) \, $.
\end{proposition}

 \Proof.
%\footnote{\bf Questa era rimasta in sospeso, c'\`e da
%meditare in particolare sulle citazioni.}
%
%
  In the classical setting, the analogous result is proved by
this argument: one uses the so-called ``straightening relations''
%(also called ``Pl{\"u}cker relations'')
to get rid of the extra
minors (see, for example,  \cite{dep},  \S 2).  Here the argument works
essentially the same, using the  {\sl quantum straightening\/
{\rm (or} Pl{\"u}cker\/{\rm )}  relations\/}  (see  \cite{fi1},  \S 4,
\cite{tt}, formula (3.2)(c) and Note I, Note II).
%the lack of commutativity will not affect the argument.
\qed

\medskip

\begin{remark}  As before, we have that
  $$  \cO_q^{\,\text{\it loc}}\big( G \big/ P \big)
\,\;{\textstyle \bigcap}\;\, (q-1) \, \cO_q^{\,\text{\it loc}}(G) \,\;
= \;\, (q-1) \, \cO_q^{\,\text{\it loc}}\big( G \big/ P \big)  $$
This can be easily deduced from Proposition  \ref{firstintersection},
taking into account Proposition  \ref{classical}.  As a consequence,
the map
  $$  \cO_q^{\,\text{\it loc}}\big( G \big/ P \big) \Big/ (q-1)
\, \cO_q^{\,\text{\it loc}}\big( G \big/ P \big) \;
\relbar\joinrel\relbar\joinrel\relbar\joinrel\longrightarrow
\; \cO_q^{\,\text{\it loc}}(G) \Big/ (q-1) \,
\cO_q^{\,\text{\it loc}}(G)  $$
is  {\sl injective},  so that the specialization map
  $$  \pi_{G/P}^{\,\text{\it loc}} \colon \,
\cO_q^{\,\text{\it loc}}\big( G \big/ P \big) \;
\relbar\joinrel\relbar\joinrel\relbar\joinrel\twoheadrightarrow
\; \cO_q^{\,\text{\it loc}}\big( G \big/ P \big) \Big/
(q-1) \, \cO_q^{\,\text{\it loc}}\big( G \big/ P \big)  $$
coincides with the restriction of the specialization map
  $$  \pi_G^{\,\text{\it loc}} \; \colon \; \cO_q^{\,\text{\it loc}}(G)
\; \relbar\joinrel\relbar\joinrel\relbar\joinrel\twoheadrightarrow
\; \cO_q^{\,\text{\it loc}}(G) \Big/ (q-1) \,
\cO_q^{\,\text{\it loc}}(G)  \quad .  $$
\end{remark}

\medskip

   The following proposition gives a description of the algebra
$ \cO_q^{\,\text{\it loc}}\big( G \big/ P \big) \, $:

\medskip

\begin{proposition} \label{bigcell}
  The big cell ring is isomorphic to a matrix algebra
%
% , namely
%
  $$  \begin{array}{cccl}
   \qquad  \cO_q^{\,\text{\it loc}}\big( G\big/P \big)  &  \lra
&  \cO_q\big(M_{(n-r) \times r}\big)  &  \\
   \qquad  t_{ij}  &  \mapsto  &  x_{ij}  &  \qquad  \forall \;\;
1 \leq j \leq r < i \leq n
      \end{array}  $$
i.e.~the generators  $ t_{ij} $'s  satisfy the Manin relations.
\end{proposition}

\Proof.  See \cite{fi2},  Proposition 1.9.   \qed

\bigskip
\medskip

\section{The Quantum Duality Principle for quantum Grassmannians}
\label{qdpgrass}

\medskip

   {} \indent   The quantum duality principle (QDP), originally due
to Drinfeld  \cite{dr}  and later formalized in  \cite{ga2}  and
extended in  \cite{ga3,ga4}  by Gavarini, is a functorial recipe
to obtain a quantum group starting from a given one.  The main
ingredients are the ``Drinfeld functors'', which are equivalences
between the category of QFA's and the category of QUEA's.  Ciccoli
and Gavarini extended this principle to the setting of homogeneous
spaces.  More precisely, in  \cite{cg}  they developed the QDP for
homogeneous spaces in the  {\sl local setting},  i.e.~for quantum
groups of formal type (where topological Hopf algebras are taken
into account).
%
%   The  {\sl global\/}  version instead is constructed
% in  \cite{cg2}.
%
 If one tries to find a global version of the QDP for
non quasi-affine homogeneous spaces, then problems arise from the very
beginning, as explained in \S \ref{intro}.
 The case of  {\sl projective\/}  homogeneous spaces has been solved in
\cite{cfg},  where the original version of the Drinfeld-like functor for
which the (global) QDP recipe should fail is suitably modified.

\medskip

   In this section, we apply the general recipe for projective homogeneous
spaces to the Grassmannian  $ G/P \, $.  The result is a quantization of
the homogeneous space  {\sl dual\/}  (in the sense of Poisson duality, see
\cite{cg})  to  $ G/P \, $,  just as the QDP recipe predicts in the setting
of  \cite{cg}.

\medskip

   We begin recalling the Drinfeld functor  $ \; {}^\vee \,
\colon \, QFA \longrightarrow QUEA \; $.

\medskip

\begin{definition}
  Let  $ G $  be an affine algebraic group over  $ \Bbbk \, $,  and
$ \cO_q(G) $  a quantization of its function algebra.  Let  $ J $
be the augmentation ideal of  $ \cO_q(G) \, $,  i.e.~the kernel
of the counit  $ \; \epsilon : \cO_q(G) \lra \Bbbk \; $.  Define
  $$  \cO_q(G)^\vee \, := \, \big\langle {(q-1)}^{-1} \, J \, \big\rangle
\, = \, {\textstyle \sum\limits_{n=0}^\infty} \, {(q-1)}^{-n} \, J^n
\qquad  \big( \subset \, \cO_q(G) \otimes_{\Bbbk_{q}} \Bbbk(q) \,\big)
\; .  $$
\end{definition}

   It turns out that  $ \cO_q(G)^\vee $  is a quantization of
$ U(\fg^*) \, $,  where  $ \fg^* $  is the dual Lie bialgebra
to the Lie bialgebra  $ \, \fg = \text{\it Lie}\,(G) \, $.  So
$ \cO_q(G)^\vee $  is a QUEA, and an infinitesimal quantization for
any Poisson group  $ G^* $  dual to  $ G \, $,  i.e.~such that  $ \,
\text{\it Lie}\,\big(G^*\big) \cong \fg^* \, $  as Lie bialgebras.
Moreover, the association  $ \, \cO_q(G) \mapsto \cO_q(G)^\vee \, $
yields a functor from QFA's to QUEA's (see  \cite{ga3,ga4}  for
more details).
%
%   $$  \begin{array}{ccc}
%       QFA  &  \lra  &  QUEA  \\
%   \cO_q(G) &  \mapsto  &  \cO_q(G)^\vee
%       \end{array}  $$
% %
% is functorial (see  \cite{ga3,ga4}  for more details).
%

\medskip

\begin{remark} \label{qdp}
  Let  $ \, G = {GL}_n \, $.  Then  $ \cO_q(G)^\vee $  is generated,
as a unital subalgebra of  $ \, \cO_q(G) \otimes_{\Bbbk_q} \Bbbk(q) \, $,
by the elements
  $$  \D_- \, := \, {(q-1)}^{-1} \, \big( D_q^{-1} - 1 \big) \; ,
\quad \hskip5pt  \chi_{ij} \, := \, {(q-1)}^{-1} \, \big( x_{ij}
- \delta_{ij} \big)   \eqno \forall \; i, j = 1, \dots, n  $$
where the  $ x_{ij} $'s  are the generators of  $ \cO_q(G) \, $.  As
$ \, x_{ij} = \delta_{ij} + (q-1) \, \chi_{ij} \in \cO_q(G)^\vee \, $,
we have an obvious embedding of  $ \cO_q(G) $  into  $ \cO_q(G)^\vee
\, $.
\end{remark}

\medskip

   In the same spirit   --- mimicking the construction in  \cite{cg}
 ---   we now want to define  $ {\cO_q\big( G \big/ P \big)}^\vee \, $
when  $ G \big/ P $  is the Grassmannian.

\medskip

   Let  $ \, G = {GL}_n \, $,  and let  $ P $  be the maximal parabolic
subgroup of \S  \ref{grass}.

\medskip

\begin{definition}  Let  $ \epsilon' $  be the natural extension to
$ \cO_q^{\,\text{\it loc}}(G/P) $  of the restriction to  $ \cO_q(G/P) $
of the counit of  $ \cO_q(G) \, $,  and let  $ \, J_{G/P}^{\,\text{\it loc}}
:= \text{\it Ker}\,(\epsilon'\,) \, $.  We define (as a subset of  $ \,
\cO_q^{\,\text{\it loc}}\big( G \big/ P \big) \otimes_{\Bbbk_q}
\Bbbk(q) \, $)
  $$  \cO_q\big( G \big/ P \big)^\vee  \, := \,  \big\langle \,
{(q-1)}^{-1} \, J_{G/P}^{\,\text{\it loc}} \, \big\rangle \, = \,
{\textstyle \sum\limits_{n=0}^\infty} \, {(q-1)}^{-n} \, {\big(
J_{G/P}^{\,\text{\it loc}} \big)}^n  \;\;\; .  $$
\end{definition}

\medskip

   It is worth pointing out that  $ \, \cO_q\big( G\big/P \big)^\vee \, $
is  {\sl not\/}  a ``quantum homogeneous space'' for  $ \cO_q(G)^\vee $
in any natural way, i.e.~it does not admit a coaction of  $ \cO_q(G)^\vee
\, $.  This is a consequence of the fact that there is no natural coaction
of  $ \cO_q(G) $  on  $ \cO_q^{\,\text{\it loc}}\big( G\big/P \big) \, $.
Now we examine this more closely.

\medskip

   Since  $ {\cO_q\big( G\big/P \big)}^\vee $  is not contained in
$ {\cO_q(G)}^\vee \, $,  we cannot have a  $ {\cO_q(G)}^\vee $  coaction
induced by the coproduct.  This would be the case if  $ {\cO_q\big(
G \big/ P \big)}^\vee $  were a (one-sided)  {\sl coideal\/}  of
$ {\cO_q(G)}^\vee \, $;  but this is not true because  $ \cO_q^{\,\text{\it
loc}}\big( G \big/ P \big) $  is not a (right) coideal of  $ \cO_q(G) $.
This reflects the geometrical fact that the big cell of  $ G/P $
is not a  $ G $--space  itself.  Nevertheless, we shall find a way
around this problem simply by {\sl enlarging}  $ {\cO_q\big( G\big/P
\big)}^\vee $  and  $ {\cO_q(G)}^\vee $,  i.e.~by taking their
$ (q \! - \! 1) $--adic  completion (which will not affect their
behavior at  $ \, q = 1 \, $).

\medskip

   To begin, we provide a concrete description of
$ {\cO_q\big( G\big/P \big)}^\vee \, $:
%
% ,  in terms of a presentation by generators and relations:
%

\medskip

\begin{proposition}
  $$  {\cO_q\big( G\big/P \big)}^\vee  \; = \;
\Bbbk_q \big\langle {\{\, \mu_{i{}j} \,\}}_{i = r+1,
\dots, n}^{j = 1, \dots, r} \big\rangle \Big/ I_M  $$
where  $ \; \mu_{i{}j} := {(q-1)}^{-1} \, t_{i{}j} \; $  (for all
$ i $  and  $ j \, $),  $ \, I_M \, $  is the ideal of the Manin
relations among the  $ \mu_{i{}j} $'s, and  $ \; t_{ij} =
{(-q)}^{r-j} \, \Dij \, D_0^{-1} \, $  (for all  $ i $
and  $ j $).
\end{proposition}

\Proof.  Trivial from definitions and Proposition \ref{bigcell}.   \qed

\bigskip

   We now explain the relation between  $ {\cO_q\big( G\big/P \big)}^\vee $
and  $ {\cO_q(G)}^\vee \, $.  The starting point is the following special
property:

\medskip

\begin{proposition}  \label{intersection}
  $$  {\cO_q\big( G \big/ P \big)}^\vee \,\;{\textstyle \bigcap}\;\,
(q-1) \, {\cO_q(G)}^\vee\big[D_0^{-1}\big]  \,\; = \;\,
(q-1) \, {\cO_q\big( G \big/ P \big)}^\vee  $$
\end{proposition}

\Proof. It is the same as for Proposition  \ref{firstintersection}.
\qed

\medskip

\begin{remark}
   As a direct consequence of Proposition  \ref{intersection},
the canonical map
   $$  {\cO_q\big( G \big/ P \big)}^\vee \! \Big/ (q-1) \,
{\cO_q\big( G\big/P \big)}^\vee \,
\relbar\joinrel\relbar\joinrel\relbar\joinrel\longrightarrow
\; {\cO_q(G)}^\vee\big[D_0^{-1}\big] \Big/ (q-1) \,
{\cO_q(G)}^\vee\big[D_0^{-1}\big]  $$
is in fact  {\sl injective\/}:  therefore, the specialization map
  $$  \pi^\vee_{G/P} \; \colon \; {\cO_q\big( G \big/ P \big)}^\vee
\; \relbar\joinrel\relbar\joinrel\relbar\joinrel\twoheadrightarrow \;
{\cO_q\big( G \big/ P \big)}^\vee \Big/ (q-1) \, {\cO_q\big(
G\big/P \big)}^\vee  $$
coincides with the restriction to  $ {\cO_q\big( G\big/P \big)}^\vee $
of the specialization map
  $$  \pi^\vee_G \; \colon \; {\cO_q(G)}^\vee\big[ D_0^{-1}\big]
\; \relbar\joinrel\relbar\joinrel\relbar\joinrel\twoheadrightarrow
\; {\cO_q(G)}^\vee\big[D_0^{-1}\big] \Big/ (q-1) \,
{\cO_q(G)}^\vee\big[D_0^{-1}\big]  \quad .  $$
\end{remark}

\medskip

   From now on, let  $ \widehat{A} $  denote the  $ (q-1) $--adic
completion of any  $ \Bbbk_q $--algebra  $ A \, $.  Note that
$ \widehat{A} $  and  $ A $  have the same specialization at
$ \, q = 1 \, $,  i.e.~$ \; A / (q-1) \, A \; $  and  $ \;
\widehat{A} / (q-1) \, \widehat{A} \; $  are canonically
isomorphic.  When  $ \, A = \cO_q(G) \, $,  \, note also that
$ \, \widehat{\cO_q(G)} \, $  is naturally a complete topological
Hopf  $ \Bbbk_q $--algebra.
                                               \par
   The next result show why it is relevant to introduce such
completions.

\medskip

\begin{lemma}  \label{embedding}
$ {\cO_q(G)}^\vee\big[D_0^{-1}\big] \, $
naturally embeds into  $ \, \widehat{{\cO_q(G)}^\vee} \, $.
\end{lemma}

\Proof.  By remark  \ref{qdp}  we have that  $ {\cO_q(G)}^\vee $  is
generated by the elements
  $$  \D_- \, := \, {(q-1)}^{-1} \, \big( D_q^{-1} - 1 \big) \; ,
\quad \hskip5pt  \chi_{ij} \, := \, {(q-1)}^{-1} \, \big( x_{ij}
- \delta_{ij} \big)   \eqno \forall \; i, j = 1, \dots, n  $$
inside  $ \, \cO_q(G) \otimes_{\Bbbk_q} \Bbbk(q) \, $.
On the other hand, observe that
  $$  \begin{array}{c}
   x_{ij} \, = \, (q-1) \, \chi_{i,j} \, \in \, (q-1) \, {\cO_q(G)}^\vee \,
\qquad  \forall \;\; i \not= j  \\
                                \\
   x_{\ell \ell} \, = \, 1 + (q-1) \, \chi_{\ell \ell} \, \in \,
\big( 1 + (q-1) \, {\cO_q(G)}^\vee \,\big) \, \qquad  \forall \;\; \ell \; .
      \end{array}   \leqno \text{and}  $$
Then, if we expand explicitly the  $ q $--determinant  $ D_0 :=
D^{I_0} \, $,  we immediately see that  $ \; D_0 \in \big( 1 +
(q-1) \, {\cO_q(G)}^\vee \,\big) \; $  as well.  Therefore  $ D_0 $  is
invertible in  $ \widehat{{\cO_q(G)}^\vee} $,  and so the natural immersion
$ \; {\cO_q(G)}^\vee \lhook\joinrel\relbar\joinrel\longrightarrow
\widehat{{\cO_q(G)}^\vee} \; $  can be canonically extended to an
immersion  $ \; {\cO_q(G)}^\vee\big[D_0^{-1}\big]
\lhook\joinrel\relbar\joinrel\longrightarrow
\widehat{{\cO_q(G)}^\vee} \; $,  \, q.e.d.   \qed

\bigskip

\begin{corollary}
\label{specializations}
                                         \hfill\break
   \indent   (a) \, The specializations at  $ \, q \! = \! 1 \, $  of\/
$ \, {\cO_q(G)}^\vee \, $,  $ \, {\cO_q(G)}^\vee\big[D_0^{-1}\big]
\, $  and  $ \, \widehat{{\cO_q(G)}^\vee} \, $  are canonically
isomorphic.  More precisely, the chain
  $$  {\cO_q(G)}^\vee \; \lhook\joinrel\relbar\joinrel\longrightarrow
\; {\cO_q(G)}^\vee\big[D_0^{-1}\big] \;
\lhook\joinrel\relbar\joinrel\longrightarrow
\; \widehat{{\cO_q(G)}^\vee}  $$
of canonical embeddings induces at  $ \, q = 1 \, $  a chain
of isomorphisms.
                                         \hfill\break
   \indent   (b) \,  $ {\cO_q\big( G \big/ P \big)}^\vee \, $  embeds
into  $ \, \widehat{{\cO_q(G)}^\vee} \, $  via the chain of embeddings
  $$  {\cO_q\big( G \big/ P \big)}^\vee \;
\lhook\joinrel\relbar\joinrel\longrightarrow
\; {\cO_q(G)}^\vee\big[D_0^{-1}\big] \;
\lhook\joinrel\relbar\joinrel\longrightarrow
\; \widehat{{\cO_q(G)}^\vee}  $$
   \indent   (c) \qquad  $ \displaystyle{
{\cO_q\big( G \big/ P \big)}^\vee \,\;{\textstyle \bigcap}\;\,
(q-1) \, \widehat{{\cO_q(G)}^\vee}  \,\; = \;\, (q-1) \,
{\cO_q\big( G \big/ P \big)}^\vee } \quad $.
\end{corollary}

\Proof.  Part  {\it (a)\/}  and  {\it (b)\/}  are trivial,
and  {\it (c)\/}  follows easily from them.   \qed

\bigskip

   Notice that part  {\it (c)\/}  of Corollary  \ref{specializations}
also implies that
  $$  {\cO_q\big( G \big/ P \big)}^\vee{\Big|}_{q=1} \,
:= \; {\cO_q\big( G \big/ P \big)}^\vee \! \Big/ (q-1)
\, {\cO_q\big( G \big/ P \big)}^\vee  $$
is a subalgebra of
  $$  \widehat{{\cO_q(G)}^\vee}{\Big|}_{q=1} \, = \;
{\cO_q(G)}^\vee{\Big|}_{q=1} \, := \; {\cO_q(G)}^\vee \Big/
(q-1) \, {\cO_q(G)}^\vee \, \cong \; U(\fg^*)  $$
just because the specialization map
  $$  \pi_{G/P}^\vee \; \colon \; {\cO_q\big( G\big/P \big)}^\vee
\; \relbar\joinrel\relbar\joinrel\relbar\joinrel\twoheadrightarrow \;
{\cO_q\big( G\big/P \big)}^\vee \Big/ (q-1) \, {\cO_q\big( G\big/P
\big)}^\vee  $$
coincides with the restriction to  $ {\cO_q\big( G\big/P \big)}^\vee $
of the specialization map
  $$  \widehat{\pi_G^\vee} \; \colon \; \widehat{{\cO_q(G)}^\vee} \;
\relbar\joinrel\relbar\joinrel\relbar\joinrel\twoheadrightarrow
\; \widehat{{\cO_q(G)}^\vee} \Big/ (q-1) \, \widehat{{\cO_q(G)}^\vee}
\quad .  $$
Now we want to see what is  $ \, {\cO_q\big( G\big/P \big)}^\vee
{\Big|}_{q=1} $  inside  $ \, U\big({\fgl_n}^{\!*}\big) \, $.
In other words, we want to understand what is the space that
$ {\cO_q\big( G\big/P \big)}^\vee $  is quantizing.

\medskip

\begin{proposition}  \label{quantumspec}
  $$  {\cO_q\big( G \big/ P \big)}^\vee{\Big|}_{q=1}
\, = \; U\big(\fp^\perp\big)  $$
as a subalgebra of  $ \; {\cO_q(G)}^\vee{\Big|}_{q=1} \! =
U\big({\fgl_n}^{\!*}\big) \, $,  where  $ \, \fp^\perp \, $
is the orthogonal subspace to  $ \, \fp := \text{\it Lie}\,(P)
\, $  inside  $ \, {\fgl_n}^{\!*} \, $.
\end{proposition}

\Proof.  Thanks to the previous discussion, it is enough to
show that
  $$  \pi_G^\vee\Big({\cO_q\big( G \big/ P \big)}^\vee\Big) \, =
\, U\big(\fp^\perp\big) \, \subseteq \, U \big( {\fgl_n}^{\!*}
\big) \, = \, {\cO_q(G)}^\vee{\Big|}_{q=1}  \quad .  $$
To do this, we describe the isomorphism  $ \, {\cO_q(G)}^\vee
{\Big|}_{q=1} \! \cong U\big({\fgl_n}^{\!*}\big) \, $  (cf.~\cite{ga1}).
First, recall that  $ {\cO_q(G)}^\vee $  is generated by the elements
(see  Remark \ref{qdp})
  $$  \D_- \, := \, {(q-1)}^{-1} \, \big( D_q^{-1} - 1 \big) \; ,
\quad \hskip5pt  \chi_{ij} \, := \, {(q-1)}^{-1} \, \big( x_{ij}
- \delta_{ij} \big)   \eqno \forall \; i, j = 1, \dots, n  $$
inside  $ \, \cO_q(G) \otimes_{\Bbbk_q} \Bbbk(q) \, $.
In terms of these  generators, the isomorphism reads
  $$  \displaylines{
   {\cO_q(G)}^\vee{\Big|}_{q=1} \,
\relbar\joinrel\relbar\joinrel\relbar\joinrel\longrightarrow
\; U\big({\fgl_n}^{\!*}\big)  \cr
   \overline{{\D}_-} \mapsto -(\e_{1,1} + \cdots + \e_{n,n})
\; ,  \qquad  \overline{\chi_{i,j}} \mapsto \e_{i,j} \qquad
\forall \;\; i \, , j \; .  \cr }  $$
where we used notation  $ \; \overline{X} := X \mod (q-1) \,
{\cO_q(G)}^\vee \; $.  Indeed, from  $ \, \overline{\chi_{i,j}}
\mapsto \e_{i,j} \, $  and  $ \, {(q-1)}^{-1} \, \big( D_q - 1
\big) \in {\cO_q(G)}^\vee \, $,  one gets  $ \, \overline{D_q}
\mapsto 1 \, $  and  $ \, \overline{{(q-1)}^{-1} \, \big( D_q
- 1 \big)} \mapsto \e_{1,1} + \cdots + \e_{n,n} \, $.  Moreover,
the relation  $ \, D_q \, D_q^{-1} = 1 \, $  in  $ \cO_q(G) $
implies  $ \; D_q \, \D_- = - {(q-1)}^{-1} \, \big( D_q - 1 \big)
\; $  in  $ {\cO_q(G)}^\vee \, $,  whence clearly  $ \; \overline{\D_-}
\mapsto -(\e_{1,1} + \cdots + \e_{n,n}) \; $  as claimed (cf.~\cite{ga1},
\S 3, or  \cite{ga3},  \S 7).

\medskip

   In other words, the specialization  $ \, \pi_G^\vee \, \colon
{\cO_q(G)}^\vee \relbar\joinrel\relbar\joinrel\twoheadrightarrow
U\big({\fgl_n}^{\!*}\big) \, $  is given by
  $$  \pi_G^\vee \big( \D_- \big) \, = \, -(\e_{1,1} + \cdots + \e_{n,n})
\; ,  \qquad  \pi_G^\vee \big( \chi_{i,j} \big) \, = \, \e_{i,j}  \qquad
\forall \;\; i \, , j \; .  $$
   \indent   If we look at  $ \widehat{{\cO_q(G)}^\vee} $,  things are
even simpler.  Since
  $$  D_q \, \in \, \Big( 1 + (q-1) \, {\cO_q(G)}^\vee \Big) \subset
\Big( 1 + (q-1) \, \widehat{{\cO_q(G)}^\vee} \Big) \; ,  $$
then  $ \, D_q^{-1} \in \Big( 1 + (q-1) \, \widehat{{\cO_q(G)}^\vee}
\Big) \, $, and the generator  $ \, \D_- \, $  can be dropped.
The specialization map  $ \widehat{\pi_{G/P}^\vee} $  of course is
still described by  formul\ae{}  as above.

\medskip

   Now let's compute  $ \, \pi_{G/P}^\vee \Big( {\cO_q\big( G\big/P
\big)}^\vee \Big) = \widehat{\pi_G^\vee} \Big( {\cO_q\big( G\big/P
\big)}^\vee \Big) \, $.  Recall that  $ {\cO_q\big( G\big/P
\big)}^\vee $  is generated by the  $ \mu_{i{}j} $'s,  with
  $$  \mu_{i{}j} \; := \; {(q-1)}^{-1} \, t_{i{}j} \; = \;
{(q-1)}^{-1} \, {(-q)}^{r-j} \, \Dij \, D_0^{-1}  $$
for  $ \, i = r+1, \dots, n \, $,  and  $ \, j = 1, \dots, r \, $;
thus we must compute  $ \, \widehat{\pi_G^\vee} \big( \mu_{i{}j}
\big) \, $.
                                           \par
   By definition, for every  $ \, i \not= j \, $  the element $ \, x_{i{}j} = (q-1) \, \chi_{i{}j} \, $  is mapped to 0 by
$ \widehat{\pi_G^\vee} \, $.  Instead, for each  $ \ell $  the
element  $ \, x_{\ell\,\ell} = 1 + (q-1) \, \chi_{\ell\,\ell} \, $
is mapped to 1 (by  $ \widehat{\pi_G^\vee} $  again).  But then,
expanding the  $ q $--determinants  one easily finds   --- much
like in the proof of Lemma  \ref{embedding}  ---   that
  $$  \displaylines{
   \widehat{\pi_G^\vee} \, \Big( {(q-1)}^{-1} \, \Dij \Big) \; = \;
\Big( {(q-1)}^{-1} \, {\textstyle \sum_{\sigma \in \, \mathcal{S}_r}}
\, {(-q)}^{\ell(\sigma)} \, x_{1 \, \sigma(1)} \cdots x_{r \, \sigma(r)}
\Big)  \, =   \hfill  \cr
   \hfill   = \; \widehat{\pi_G^\vee} \, \Big( {(q-1)}^{-1} \, {\textstyle
\sum\limits_{\sigma \in \, \mathcal{S}_r}} {(-q)}^{\ell(\sigma)} \,
\big( \delta_{1 \, \sigma(1)} + (q-1) \chi_{1 \, \sigma(1)}) \cdots
\big( \delta_{1 \, \sigma(r)} + (q-1) \, \chi_{1 \, \sigma(r)})
\Big)  \cr }  $$
The only term in  $ (q-1) $  in the expansion of  $ \Dij $  comes
from the product
  $$  \big( 1 + (q \! - \! 1) \chi_{1 \, 1}) \cdots \big( 1 +
(q \! - \! 1) \chi_{r \,r} \big) \, (q \! - \! 1) \chi_{i \, j}
\equiv  (q \! - \! 1) \chi_{i \, j}  \! \mod {(q \! - \! 1)}^2
\cO\big(G \! \big/ \! P\big)  $$
Therefore, from the previous analysis we get
  $$  \displaylines{
   \widehat{\pi_G^\vee} \, \Big( {(q-1)}^{-1} \, \Dij \Big)
\; = \; \widehat{\pi_G^\vee} \, \big( \chi_{i,j} \big) \;
= \; \e_{i,j}  \cr
   \widehat{\pi_G^\vee} \, \big( D_0 \big) \; = \;
\widehat{\pi_G^\vee} \big( 1 \big) \; = \; 1 \;\; ,  \qquad
\widehat{\pi_G^\vee} \, \big( D_0^{-1} \big) \; = \;
\widehat{\pi_G^\vee} \big(1 \big) \; = \; 1  \cr }  $$
hence we conclude that  $ \; \widehat{\pi_G^\vee} \big( \mu_{i{}j} \big)
= {(-1)}^{r-j} \, \e_{i,j} \; $,  \, for all  $ \, 1 \leq j \leq r < i
\leq n \, $.

\medskip

   The outcome is that  $ \; \pi_{G/P}^\vee \Big( {\cO_q\big( G\big/P
\big)}^\vee \Big) \, = \, U(\mathfrak{h}) \; $,  \; where
  $$  \mathfrak{h} \, := \, \text{\it Span}\,\big( \big\{\, \e_{i,j}
\,\big|\, r+1 \leq i \leq n \, , \; 1 \leq j \leq r \,\big\} \big)
\;\; .  $$
On the other hand, from the very definitions and our description
of  $ \, {\fgl_n}^{\!*} \, $  one easily finds that
$ \; \mathfrak{h} = \fp^\perp \; $,  \; for  $ \; \fp :=
\text{\it Lie}\,(P) \; $.  The claim follows.   \qed

\bigskip

   Proposition  \ref{quantumspec}  claims that  $ {\cO_q \big(
G\big/P \big)}^\vee $  is a quantization of  $ U \big( \fp^\perp
\big) \, $,  i.e.~it is a unital  $ \Bbbk_q $--algebra  whose
semiclassical limit is  $ U\big(\fp^\perp\big) \, $.  Now, the
fact that  $ U\big(\fp^\perp\big) $  describes (infinitesimally)
a homogeneous space for  $ G^* $  is encoded in algebraic terms
by the fact that it is a (left) coideal of  $ U(\fg^*) \, $;
in other words,  $ U\big(\fp^\perp\big) $  is a (left)
$ U(\fg^*) $--comodule  w.r.t.~the restriction of the coproduct
of  $ U(\fg^*) \, $.  Thus, for  $ {\cO_q \big( G\big/P \big)}^\vee $
to be a quantization of  $ U\big(\fp^\perp\big) $  {\sl as a
homogeneous space\/}  we need also a quantization of this fact:
namely, we would like  $ {\cO_q \big( G\big/P \big)}^\vee $  to
be a left coideal of  $ {\cO_q(G)}^\vee $,  our quantization of
$ U(\fg^*) \, $.  But this makes no sense at all,
       \hbox{as  $ {\cO_q \big(
G \! \big/ \! P \big)}^\vee $  is not even a subset of
$ {\cO_q (G)}^\vee \, $!}
                                            \par
   This problem leads us to enlarge a bit our quantizations  $ {\cO_q
\big( G\big/P \big)}^\vee $  and  $ {\cO_q(G)}^\vee \, $:  we take
their  $ (q \! - \! 1) $--adic  completions, namely  $ \widehat{{\cO_q
\big( G\big/P \big)}^\vee} $  and  $ \widehat{{\cO_q(G)}^\vee} \, $.
While not affecting their behavior at  $ \, q = 1 \, $  (i.e., their
semiclassical limits are the same), this operation solves the problem.
Indeed,  $ \widehat{{\cO_q(G)}^\vee} $  is big enough to contain
$ {\cO_q \big( G\big/P \big)}^\vee \, $,  by Corollary
\ref{specializations}{\it (b)}.  Then, as  $ \widehat{{\cO_q(G)}^\vee} $
is a topological Hopf algebra, inside it we must look at the closure
of  $ {\cO_q\big( G\big/P \big)}^\vee \, $.  Thanks to Corollary
\ref{specializations}{\it (c)\/}  (which means, roughly, that an
Artin-Rees lemma holds), the latter is nothing but  $ \widehat{{\cO_q
\big( G\big/P \big)}^\vee} $.  Finally, next result tells us that
$ \widehat{{\cO_q \big( G\big/P \big)}^\vee} $  is a left coideal
of  $ \widehat{{\cO_q(G)}^\vee} $,  as expected.

\bigskip

\begin{proposition}  $ \, \widehat{{\cO_q\big( G \big/ P \big)}^\vee} \, $
is a left coideal of  $ \; \widehat{{\cO_q(G)}^\vee} \, $.
\end{proposition}

\Proof.  Recall that the coproduct  $ \widehat{\Delta} $  of
$ \widehat{{\cO_q(G)}^\vee} $  takes values in the  {\sl topological\/}
tensor product  $ \, \widehat{{\cO_q(G)}^\vee} \,\widehat{\otimes}\,
\widehat{{\cO_q(G)}^\vee} \, $,  which by definition is the  $ (q \!
- \! 1) $--adic  completion of the  {\sl algebraic\/}  tensor product
$ \, \widehat{{\cO_q(G)}^\vee} \otimes \widehat{{\cO_q(G)}^\vee} \, $.
Our purpose then is to show that this coproduct  $ \widehat{\Delta} $
maps  $ \widehat{{\cO_q \big( G\big/P \big)}^\vee} $  in the topological
tensor product  $ \, \widehat{{\cO_q(G)}^\vee} \,\widehat{\otimes}\,
\widehat{{\cO_q \big( G\big/P \big)}^\vee} \, $.
                                          \par
   By construction, the coproduct of  $ {\cO_q(G)}^\vee $,  hence of
$ \widehat{{\cO_q(G)}^\vee} $  too, is induced by that of  $ \cO_q(G)
\, $,  say  $ \; \Delta : \cO_q(G) \lra \cO_q(G) \otimes \cO_q(G) \; $.
Now, the latter can be uniquely (canonically) extended to a
coassociative algebra morphism
  $$ \tdelta : \, \cO_q(G)\big[D_{{I_0}}^{-1}\big] \relbar\joinrel\lra
\cO_q(G)\big[D_{{I_0}}^{-1}\big] \,\widetilde{\otimes}\,
\cO_q(G)\big[D_{{I_0}}^{-1}\big]  $$
where  $ \, \widetilde{\otimes}\, $  is the  $ J_\otimes $--adic
completion of the algebraic tensor product, with
  $$  J_\otimes \; := \; J \otimes \cO_q(G) \, + \, \cO_q(G) \otimes J
\;\; ,  \qquad  J \, := \text{\it Ker}\,\big(\epsilon_{\cO_q(G)}\big)
\;\; .  $$
In fact, since  $ \; \Delta(D_0) = D_0 \otimes D_0 +
{\textstyle \sum\limits_{K \neq {I_0}}} D^{I_0}_K \otimes D^K \; $,
\, one easily computes
  $$  \displaylines{
   \tdelta(D_0)  \; = \;  \Big(\, 1 + {\textstyle \sum\limits_{K \neq
{I_0}}} \, D_K^{I_0} \, D_0^{-1} \otimes D^K D_0^{-1} \,\Big)
\big( D_0 \otimes D_0 \big)  \cr
  \tdelta\big(D_0^{-1}\big)  \; = \;  {\big( D_0 \otimes D_0 \big)}^{-1}
{\Big(\, 1 + {\textstyle \sum\limits_{K \neq {I_0}}} \, D_K^{I_0} \,
D_0^{-1} \otimes D^K D_0^{-1} \,\Big)}^{-1}  \cr
  \hfill   = \,  \Big( D_0^{-1} \otimes D_0^{-1} \Big) \;
{\textstyle \sum\limits_{n \geq 0}} \, {(-1)}^n \, {\bigg(\,
{\textstyle \sum\limits_{ K \neq {I_0}}} \, D_K^{I_0} \,
D_0^{-1} \otimes D^K D_0^{-1} \,\bigg)}^{\!n}  \cr }  $$
   \indent   Let's now look at the restriction  $ \tdelta_r $  of
$ \tdelta $  to  $ \cO_q^{\,\text{\it loc}}\big( G\big/P \big) \, $.
We have
  $$  \displaylines{
   \tdelta_r(t_{ij})  \; = \;  \tdelta_r \big( \Dij \, D_0^{-1}
\,\big)  \; = \; \tdelta\big(\Dij\big) \cdot {\tdelta\big( D_0
\big)}^{-1}  \; =   \hfill  \cr
   \hfill   = \,  \bigg(\, {\textstyle \sum\limits_L} \,
\Dij_{\,L} \, D_0^{-1} \otimes D^L \, D_0^{-1} \,\bigg) \cdot
{\textstyle \sum\limits_{n \geq 0}} \, {(-1)}^n \, {\bigg(\,
{\textstyle \sum\limits_{K \neq {I_0}}} \, D_K^{I_0} \,
D_0^{-1} \otimes D^K D_0^{-1} \,\bigg)}^{\!n}  \cr }  $$
   \indent   Now, by Proposition  \ref{classical}  we know that each
product  $ \, D^L \, D_{I_0}^{-1} \, $  is a combination of the
$ t_{ij} $'s.  Hence the formula above shows that  $ \tdelta_r $
maps  $ \cO_q^{\,\text{\it loc}}\big( G\big/P \big) $  into  $ \;
\cO_q(G)\big[D_0^{-1}\big] \,\widetilde{\otimes}\,
\cO_q^{\,\text{\it loc}}\big( G\big/P \big) \; $.
                                                 \par
   By scalar extension,  $ \tdelta $  uniquely extends to a
map defined on the  $ \Bbbk(q) $--vector  space  $ \, \Bbbk(q)
\otimes_{\Bbbk_q} \cO_q(G)\big[D_0^{-1}\big] \, $,  which
we still call  $ \tdelta \, $.  Its restriction to the similar
scalar extension of  $ \cO_q^{\,\text{\it loc}}\big( G\big/P \big) $
clearly coincides with the scalar extension of  $ \tdelta_r \, $,
hence we call it  $ \tdelta_r $  again.  Finally, the restriction of
$ \tdelta $  to  $ {\cO_q(G)}^\vee\big[D_0^{-1}\big] $  and of
$ \tdelta_r $  to  $ {\cO_q\big( G\big/P \big)}^\vee $  both coincide
--- by construction ---   with the proper restrictions of the coproduct
of  $ \widehat{{\cO_q(G)}^\vee} $  (cf.~Corollary  \ref{specializations}).
                                                 \par
   In the end, we are left to compute  $\tdelta_r(\mu_{ij}) \, $.  The
computation above gives
  $$  \displaylines{
\widehat{\Delta}(\mu_{ij})  \; = \;  \tdelta_r(\mu_{ij})  \; = \;
{(q \! - \! 1)}^{-1} \, \tdelta_r(t_{ij})  \; = \;   \hfill  \cr
   \hfill   =  {(q \! - \! 1)}^{-1} \, {\textstyle \sum\limits_L}
\, \Dij_{\,L} \, D_0^{-1} \otimes D^L \, D_0^{-1} \cdot
{\textstyle \sum\limits_{n \geq 0}} {(-1)}^n {\bigg({\textstyle
\sum\limits_{K \neq {I_0}}} \! D^{I_0}_K \, D_0^{-1} \otimes
D^K D_0^{-1} \bigg)}^{\!n}  \cr }  $$
Now, each left-hand side factor above belongs to  $ \, \widehat{{\cO_q(G)}^\vee}
\,\widehat{\otimes}\, \widehat{{\cO_q \big( G\big/P \big)}^\vee} \, $,  because
either  $ \, D^L \in J_{G/P}^{\,\text{\it loc}} \, $  (if  $ \, L \not= I_0 \, $,
with notation of \S 4.3), or  $ \, \Dij_{\,L} \in J \, $  (if  $ \, L = I_0 \, $,
with  $ \, J := \text{\it Ker}\, \big(\epsilon_{\cO_q(G)}\big) \, $).
On right-hand side instead we have
  $$  D^K \, \in \, J_{G/P}^{\,\text{\it loc}} \, \subseteq \,
(q-1) \, {\cO_q\big( G\big/P \big)}^\vee \;\; ,  \qquad  D_K^{I_0}
\, \in \, J \, \subseteq \, (q-1) \, {\cO_q(G)}^\vee  $$
whence   --- as  $ \, D_0^{-1} \in \widehat{{\cO_q(G)}^\vee} \, $  and
$ \, D_0^{-1} \in \widehat{{\cO_q\big( G\big/P \big)}^\vee} \, $  ---
we get
  $$   \! {\textstyle \sum\limits_{K \neq {I_0}}}
D_K^{I_0} \, D_0^{-1} \otimes D^K D_0^{-1}  \; \in \;
{(q-1)}^2 \, \widehat{{\cO_q(G)}^\vee} \,\widehat{\otimes}\,
\widehat{{\cO_q\big( G\big/P \big)}^\vee}  $$
so that
  $ \quad \displaystyle{
   {\textstyle \sum\limits_{n \geq 0}} \, {(-1)}^n \, {\bigg(\! {\textstyle
\sum\limits_{K \neq {I_0}}} D^{I_0}_K \, D_0^{-1} \otimes D^K D_0^{-1}
\bigg)}^{\!n}  \, \in \;  \widehat{{\cO_q(G)}^\vee} \,\widehat{\otimes}\,
\widehat{{\cO_q\big( G\big/P \big)}^\vee}  \quad . } $
                                                           \par
   The final outcome is  $ \; \widehat{\Delta}(\mu_{ij})
\in \widehat{{\cO_q(G)}^\vee} \,\widehat{\otimes}\,
\widehat{{\cO_q\big( G\big/P \big)}^\vee} \; $  for all  $ i $,  $ j \, $.
As the  $ \mu_{ij} $'s  topologically generate  $ \, \widehat{{\cO_q \big(
G\big/P \big)}^\vee} \, $,  \, this proves the claim.   \qed

\bigskip

   In the end, we get the main result of this paper.

\medskip

\begin{theorem}  \label{maintheorem}
  $ \, \widehat{{\cO_q\big( G \big/ P \big)}^\vee} \, $
is a quantum homogeneous  $ \, G^* $--space,  which is
an infinitesimal quantization of the homogeneous  $ \,
G^* $--space  $ \fp^\perp \, $.
\end{theorem}

\Proof.  Just collect the previous results.  By Proposition
\ref{quantumspec}  and by the fact that  $ \; \widehat{{\cO_q
\big( G\big/P \big)}^\vee}{\Big|}_{q=1} = \, {\cO_q\big( G\big/P
\big)}^\vee{\Big|}_{q=1} \; $ we have that the specialization of
$ \widehat{{\cO_q\big( G\big/P \big)}^\vee} $  is  $ U \big(
\fp^\perp \big) \, $.  Moreover we saw that  $ \widehat{{\cO_q
\big( G\big/P \big)}^\vee} $  is a subalgebra, and left coideal,
of  $ \widehat{{\cO_q(G)}^\vee} \, $.  Finally, we have
  $$  \widehat{{\cO_q\big( G\big/P \big)}^\vee} \,\;{\textstyle
\bigcap}\;\, (q-1) \, \widehat{{\cO_q(G)}^\vee}  \,\; = \;\, (q-1)
\, \widehat{{\cO_q\big( G\big/P \big)}^\vee}  $$
as an easy consequence of Corollary \ref{specializations}  {\it (c)}.
Therefore,  $ \widehat{{\cO_q\big( G\big/P \big)}^\vee} $  is a
quantum homogeneous space, in the usual sense.
As  $ \widehat{{\cO_q(G)}^\vee} $  is a quantization of  $ \fg^* \, $,
we have that  $ \widehat{{\cO_q\big( G \big/ P \big)}^\vee} $  is
in fact a quantum homogeneous space for  $ G^* \, $;  of course, this
is a quantization of  {\sl infinitesimal\/}  type.
%
% because we are quantizing a
% universal enveloping algebra and subobjects
% (subalgebras, coideals) of it.
%
 \qed

\bigskip

\begin{remark}
  All these computations can be repeated, step by step, taking
$ \, G = {SL}_n \, $  and  $ \, P = P_1 \, $.
\end{remark}

\end{document}